\documentclass[11pt,a4paper,reqno]{amsart} 

\usepackage{amssymb,graphicx}

\usepackage{amssymb}

\newcounter{mnotecount}[section]

\renewcommand{\themnotecount}{\thesection.\arabic{mnotecount}}

\newcommand{\mnote}[1]
{\protect{\stepcounter{mnotecount}}$^{\mbox{\footnotesize
$
\bullet$\themnotecount}}$ \marginpar{
\raggedright\tiny\em
$\!\!\!\!\!\!\,\bullet$\themnotecount: #1} }

\newcommand{\hook}{{\setlength{\unitlength}{11pt}   
                   \begin{picture}(.833,.8)
                   \put(.15,.08){\line(1,0){.35}}
                   \put(.5,.08){\line(0,1){.5}}
                   \end{picture}}}
\newcommand{\CP}{\mathbb{CP}}

\newcommand{\PP}{\mathbb{P}}
\newcommand{\RP}{\mathbb{RP}}
\newcommand{\R}{\mathbb{R}}

\def\p{\partial}
\def\be{\begin{equation}}
\def\ee{\end{equation}}

\def\bea{\begin{eqnarray}}
\def\eea{\end{eqnarray}}

\begin{document}\date{February, 2011}
\title{On the 7th order ODE with submaximal symmetry}
\author{Maciej Dunajski}
\address{Department of Applied Mathematics and Theoretical Physics\\ 
University of Cambridge\\ Wilberforce Road, Cambridge CB3 0WA\\ UK.}
\email{m.dunajski@damtp.cam.ac.uk}
\author{Vladimir Sokolov}
\address{Landau Institute for Theoretical Physics\\ Moscow\\Russia.}
\email{vsokolov@landau.ac.ru}

\begin{abstract} 
We find a general solution to the unique 7th order ODE admitting ten dimensional group of contact symmetries. The integral curves of this ODE are rational 
contact curves in $\PP^3$ which give rise to rational plane curves
of degree six. The moduli space of these curves is a real form of the homogeneous space $Sp(4)/SL(2)$.

\end{abstract}   
\maketitle
The 7th order ODE 
\be
\label{equation}
10 (y^{(3)})^3y^{(7)}- 70(y^{(3)})^2y^{(4)}y^{(6)}
-49(y^{(3)})^2(y^{(5)})^2+280(y^{(3)})(y^{(4)})^2y^{(5)}
-175 (y^{(4)})^4=0,
\ee
where $y=y(x)$ and $y^{(k)}=d^ky/dx^k$      
has recently appeared explicitly \cite{S88, olver, Doubrov2, DGN10} or implicitly \cite{Bryant1}
in several different contexts. It is the unique (up to contact transformations) equation admitting  ten--dimensional algebra of contact symmetries \cite{S88}, 
and the aim of this note is to show that its general solution is given by 
a degree six rational curve of the form
\be
y^3+\alpha(x)y^2 +\beta(x)y+\gamma(x)=0,
\label{sol_introduction}
\ee
where $(\alpha, \beta, \gamma)$ are a quadratic, a quartic, and a sextic respectively
with the coefficients depending  on seven parameters as in formula
(\ref{final_formula}).

In fact the symmetry algebra of 
(\ref{equation}) was known to Lie \cite{Lie2} who also proved that this is the maximal algebra of contact
vector fields on the plane. It is quite possible that
equation (\ref{equation}) and its general solution were also known to Lie. 
We have been unable to find it in any of Lie's works , but we 
would be grateful to hear from anyone who has (the earliest reference to the equation - but not its solution - we are aware of is \cite{Noth}).
\subsection*{Contact Lie algebras}
Let $U\subset \R^2$ be an open set and let $\PP(T^*U)$ be a projectivised 
cotangent bundle with a contact one--form $\omega$. A curve $\gamma\subset \PP(TU)$ is called {\em contact} if $\omega|_\gamma=0$.
A contact transformation is a map $f:\PP(TU)\rightarrow \PP(TU)$ which
takes contact curves into contact curves. Equivalently
$f^*(\omega)=\lambda\omega$ for some function $\lambda$. Let
$(x, y)$ be the local coordinates on $U$ and let $z$ parametrise the fibers of $\PP(TU)$, so that we can set $\omega=dy-zdx$. Consider a one parameter group of contact transformations. Close to the identity, this is characterised by a contact vector field
$X$ such that the contact condition ${\mathcal L}_X\omega=c\omega$ holds, 
where 
${\mathcal L}_X=d(X\hook)+ X\hook d$ is the Lie derivative and $c$ is some 
function. 
The contact condition implies that locally there exist a function
$H=H(x, y, z)$ such that
\be
\label{contact_vector}
X_H=-(\p_z H)\p_x+(H-z\p_z H)\p_y+(\p_x H+z\p_y H)\p_z.
\ee
If $H=a(x, y)+z b(x, y)$ then $X_H$ generates a prolongation of a family of
{\em point transformations} $f:U\rightarrow U$. Otherwise it generates a proper contact flow. 
\subsection*{The symmetry}
A remarkable result of Lie is that a maximum dimension of a Lie algebra of 
proper contact vector fields on the plane is 
ten. This maximal, ten--dimensional Lie algebra is generated by vector fields
(\ref{contact_vector}) corresponding to functions
\be
\label{10Dalgebra}
1, x, x^2, y, z, xz, x^2z-2xy, z^2, 2yz-xz^2, 4xyz-4y^2-x^2z^2.
\ee
This algebra is isomorphic to $\mathfrak{so}(5)$, or equivalently
to $\mathfrak{sp}(4)$ (this is, up to a choice of the real form, the maximal subalgebra 
of the 11 dimensional symmetry algebra of the trivial 7th order ODE). 
The first seven generators correspond to prolonged point vector fields and the last three are the proper contact vector fields.

The algebra (\ref{10Dalgebra}) generates all contact symmetries of equation (\ref{equation}): if $y(x)$ is a solution to (\ref{equation}) then so is $\tilde{y}(\tilde{x})$, where 
\[
\tilde{x}=\tilde{x}(x, y, z, c_i), \quad
\tilde{y}=\tilde{y}(x, y, z, c_i), \quad
\tilde{z}=\tilde{z}(x, y, z, c_i)=\frac{\p_x\tilde{y}+z\p_y\tilde{y}}{\p_x\tilde{x}
+z\p_y\tilde{x}}, \quad \p_z\tilde{y}=\tilde{z}\p_z\tilde{x}
\]
is the contact transformation generated by the vector fields 
(\ref{contact_vector}) and $c_1, \dots, c_{10}$ are parameters of this transformation. Equivalently, any of the generators
(\ref{10Dalgebra}) satisfies the linearisation of (\ref{equation})
when $z=y'$.

Another result of Lie \cite{Lie2} is that a maximal dimension of the
contact symmetry algebra of an ODE of order $n>3$ is $(n+4)$, with maximal symmetry occurring if only if the ODE is contact equivalent to
a trivial equation $y^{(n)}=0$. Therefore equation (\ref{equation})
is of submaximal type \cite{olver} - it is not equivalent to the 
trivial equation and 
its symmetry algebra has the largest possible dimension. Up to the contact equivalence (\ref{equation}) is the unique 7th order ODE with this property.
\subsection*{The solution}
We verify that the algebraic curve  
$
y^2+x(x-1)^3=0
$
solves (\ref{equation}). This curve is parametrised by
\be
\label{par1}
x(t)=\frac{1}{t^2+1},\quad y(t)=-\frac{t^3}{(t^2+1)^2}.
\ee
The 7--dimensional subalgebra of (\ref{10Dalgebra})  consisting of prolonged 
point symmetries integrates to 
\be
\label{point_tr}
y\rightarrow c_4y+ c_1+c_2x+c_3x^2, \quad x\rightarrow c_5 x+c_6,\quad\mbox{and}\quad
y\rightarrow \frac{y}{(1+c_7 x)^2}, \quad  x\rightarrow \frac{x}{1+c_7 x}.
\ee
The contact 
transformations generated by (\ref{contact_vector}) with
$H_8=z^2, H_9=2yz-xz^2$  and $H_{10}=4xyz-4y^2-x^2z^2$  respectively are given by
\begin{eqnarray}
\label{cont_trans_new}
\tilde{x}&=&x-2c_8 z, \quad
\tilde{y}=y-c_8 z^2, \quad \tilde{z}=z\\
\tilde{x}&=&\frac{x(1+c_9z)-2c_9y}{1-c_9z}, \quad
\tilde{y}=\frac{y(1-2c_9z)+c_9xz^2}{(1-c_9z)^2},\quad 
\tilde{z}=\frac{z}{1-c_9z},\qquad
\mbox{and}\nonumber\\
\tilde{x}&=&\frac{x}{1+4c_{10}y-2c_{10}xz},\quad
\tilde{y}=
\frac{y+4c_{10}y^2-4c_{10}xyz+c_{10}x^2z^2}{(1+4c_{10}y-2c_{10}xz)^2}, 
\, \tilde{z}=\frac{z}{1+4c_{10}y-2c_{10}xz}\nonumber
\end{eqnarray}
where $c_8, c_9, c_{10}\in \R$.

Applying the group of point transformations to the given solution  yields the 
six parameter family of solutions to (\ref{equation}) given by 
a family of algebraic curves of degree four
\be
\label{degree_six}
(y+Q)^2+P=0
\ee
where $Q=Q(x)$ is an arbitrary quadratic, and $P=P(x)$ is a quartic with one simple and one triple root. 

Recall \cite{alg_geom} that a point $p$ on 
an algebraic curve $f(x, y)=0$  in $\CP^2$ is singular if the partial 
derivatives $f_x$ and $f_y$ vanish at $p$. Moreover $p$ has multiplicity $m$
if  all  $(m-1)$st derivatives of $f$ vanish at $p$ but at least 
one $m$th derivative does not. A singular point is called ordinary if
the tangents to all branches at the point are distinct.
Any singular point characterised by a triple of 
integers $(m, \delta, r)$, where $m$ is the multiplicity, 
$r$ is the number of branches and $\delta$ is the number
of multiplicity two ordinary singular  points concentrating at $p$.
The arithmetic genus of a curve is given by
\[
{\tt g}=\frac{(d-1)(d-2)}{2}-\sum \delta,
\]
where $d$ is the degree of the curve, and the 
summation is taken over all singular points.
The curves in the family 
(\ref{degree_six})
have two singularities:
a cusp of type $(2, 1, 1)$ at $(x, y)=(x_0, -Q(x_0))$, where $x_0$ is the 
triple root of  $P$ and a point at $\infty$ of type (2, 2, 2).
Calculating the genus of curves
in this family yields
\[
{\tt g}=\frac{3\cdot 2}{2}-1-2=0.
\]
Therefore the family is rational. Applying the point 
transformations in  (\ref{point_tr}) to 
the parametrisation (\ref{par1}) and redefining the constants 
$c_1, \dots, c_6$ we find that the rational parametrisation is given by
\[
x(t)=\frac{b_5+b_6t^2}{b_0+t^2}, \quad y(t)=\frac{b_4t^4+b_3t^2+b_2 t+b_1}{(b_0+t^2)^2}.
\]
The six parameters in (\ref{degree_six}) are algebraic expressions
in the seven parameters $(b_0, \dots, b_6)$ one of which is irrelevant
and arises only in the parametrisation.

The family of curves (\ref{degree_six}) can also be obtained applying $H_9$ and the point transformations to the trivial solution $y=x^2$. It is not the general solution to (\ref{equation}) as it depends on six parameters rather than seven. To introduce the additional parameter, and construct the general
solution we use the contact transformation (\ref{cont_trans_new}) generated 
by
$(\ref{contact_vector})$ with $H_8=z^2$
\be
\label{parametrisation_curves}
\tilde{x}(t)=x(t)-2b z(t), \quad
\tilde{y}(t)=y(t)-b z(t)^2, \quad \tilde{z}(t)=z(t),
\ee
where 
\[
z(t)=\frac{\dot{y}(t)}{\dot{x}(t)}=\frac{(4b_4b_0-2b_3)t^3-3b_2t^2+
(2b_3b_0-4b_1)t+b_2b_0}{2(b_5-b_0b_6)(b_0t+t^3)},
\]
and $(b, b_0, \dots, b_6)$ are constant parameters.
The relation (\ref{parametrisation_curves}) gives a  
seven--dimensional
family of
rational contact curves in $\PP^3$. The connected component of the 
symplectic group $Sp(4)$
acts on $\PP^3$ and preserves the family
(\ref{parametrisation_curves}). The symmetry group of any fixed rational curve in this family is $SL(2)$, and so the seven--dimensional space of solutions to (\ref{equation})
is the symmetric space $Sp(4)/SL(2)$. In the holomorphic category, a rational 
curve in $\PP^3$ can be characterised by a normal bundle, which in our case is
$N={\mathcal O}(5)\oplus{\mathcal O}(5)$. The contact modification  of the 
Kodaira theorem described in  \cite{Bryant2} can be applied to deduce that we have 
constructed  a complete analytic family of contact curves: an infinitesimal contact
deformation of any fixed curve in the family (\ref{parametrisation_curves})
also belongs to this family.
\vskip5pt

Alternatively,  the general solution to
(\ref{equation}) can be given by an implicit relation
\be
\label{implicit_uv}
u(x, y, z)=0, \qquad v(x, y, z)=0
\ee
where $(x, y, z)$ are coordinates on an open set in $\PP^3$. 
Using (\ref{degree_six}) and (\ref{parametrisation_curves}) we find
\[
u=(y+bz^2+Q(X))^2+P(X) , \quad v=4P(z+Q')^2+(P')^2
\]
where $Q=Q(X)$ is a quadratic, $P=P(X)$ is a quartic with one simple
and one triple root, and $X=x+2bz$.
To find an explicit formula for $y(x)$ pick a real root $z$ 
of the cubic $v=0$, substitute this in $u=0$ and solve the 
resulting quadratic for $y$. Alternatively we use the 
resultant to produce a  planar curve birationaly equivalent to (\ref{implicit_uv}) by eliminating $z$ between $u$ and $v$. Recall \cite{Grace_Young}
that a resultant of
two polynomials
\[
u(z)=u_0+u_1z+\dots+u_m z^m, \quad v(z)=v_0+v_1z+\dots+v_n z^n
\]
is the determinant of the matrix
\[
\left( \begin{array}{cccccccc}
u_m & u_{m-1} &\cdots  & \cdots & u_0   &        &       &     \\
    &u_m      & u_{m-1}& \cdots &\cdots & u_0    &       &     \\ 
    &         &\ddots  &        &       &        &\ddots &      \\
    &         &        &  u_m   &u_{m-1}&\cdots  &\cdots & u_0 \\
v_n & v_{n-1} &\cdots  & \cdots & v_0   &        &       &     \\
    &v_n      & v_{n-1}& \cdots &\cdots & v_0    &       &     \\ 
    &         &\ddots  &        &       &        &\ddots &     \\
    &          &       &  v_n   &v_{n-1}&\cdots  &\cdots & v_0  
\end{array} \right).
\]
The resultant vanishes if $u, v$ admit  a common root.
To obtain a manageable formula consider (\ref{implicit_uv}) with
$Q=0$ and $P=X(X-1)^3$, and find the resultant of $u$ and $v$. This resultant
factorises into two terms, each giving a rational curve of degree six. 
We choose one
of these two curves (The second curve is not a solution to the 
ODE (\ref{equation}). It arises because the expressions in 
(\ref{implicit_uv}) are squares of actual solutions and  contain the term $z^2$.
Only one of the roots satisfies $z=y'$).
\begin{eqnarray}
\label{new_curve}
&&\left( 64\,b+1024\,{b}^{3} \right) {y}^{3}+ \left(  \left( 768\,{b}^{
2}+16 \right) {x}^{2}-768\,x{b}^{2}+288\,{b}^{2} \right) {y}^{2}\nonumber\\
&&+
\left( 264\,{x}^{2}b-108\,{b}^{3}+192\,{x}^{4}b-72\,xb-384\,{x}^{3}b
 \right) y\\
&&+\left(48\,{x}^{4}-27\,{b}^{2}+54\,x{b}^{2}-16\,{x}^{3}-27\,{x}^{2}{b}^{2}-48\,{x}^{5}+16\,{x}^{6}\right)=0\nonumber.
\end{eqnarray}
This curve  does not belong to the class (\ref{degree_six}):
it has two $(2, 1, 1)$ cusps and one non-ordinary $(2, 2, 2)$ singularity 
at $\infty$ whereas
$(\ref{degree_six})$ has one cusp (apart from 
the singularity at $\infty$).
\begin{figure}
\caption{Curve 
(\ref{new_curve}) with $b=1/2$ and curve 
(\ref{degree_six}) with $Q=0, P=x(x-1)^3$.}
\label{IS_levelsurface}
\begin{center}
\includegraphics[width=4cm,height=7cm,angle=270]{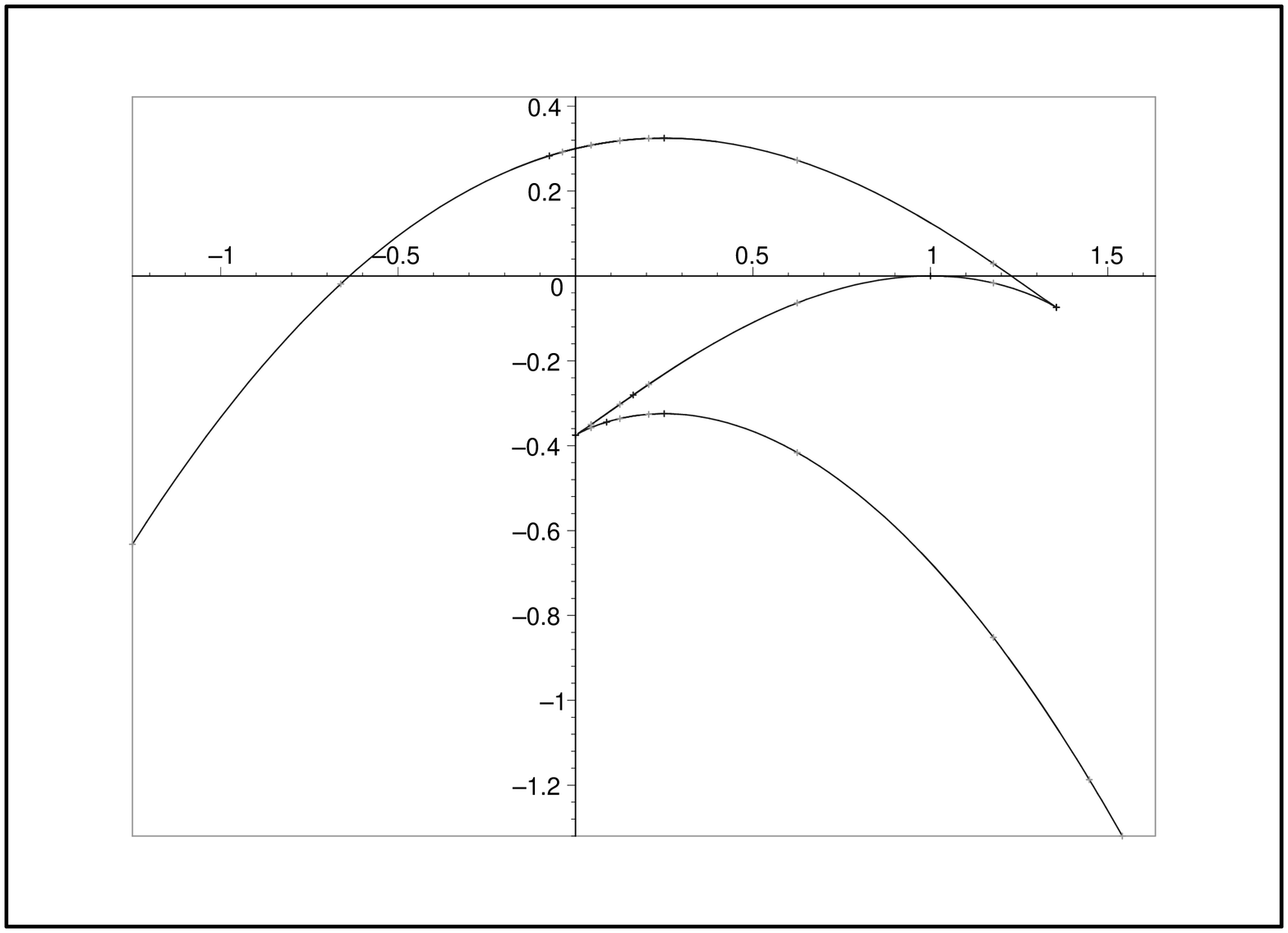},\quad
\includegraphics[width=4cm,height=7cm,angle=270]{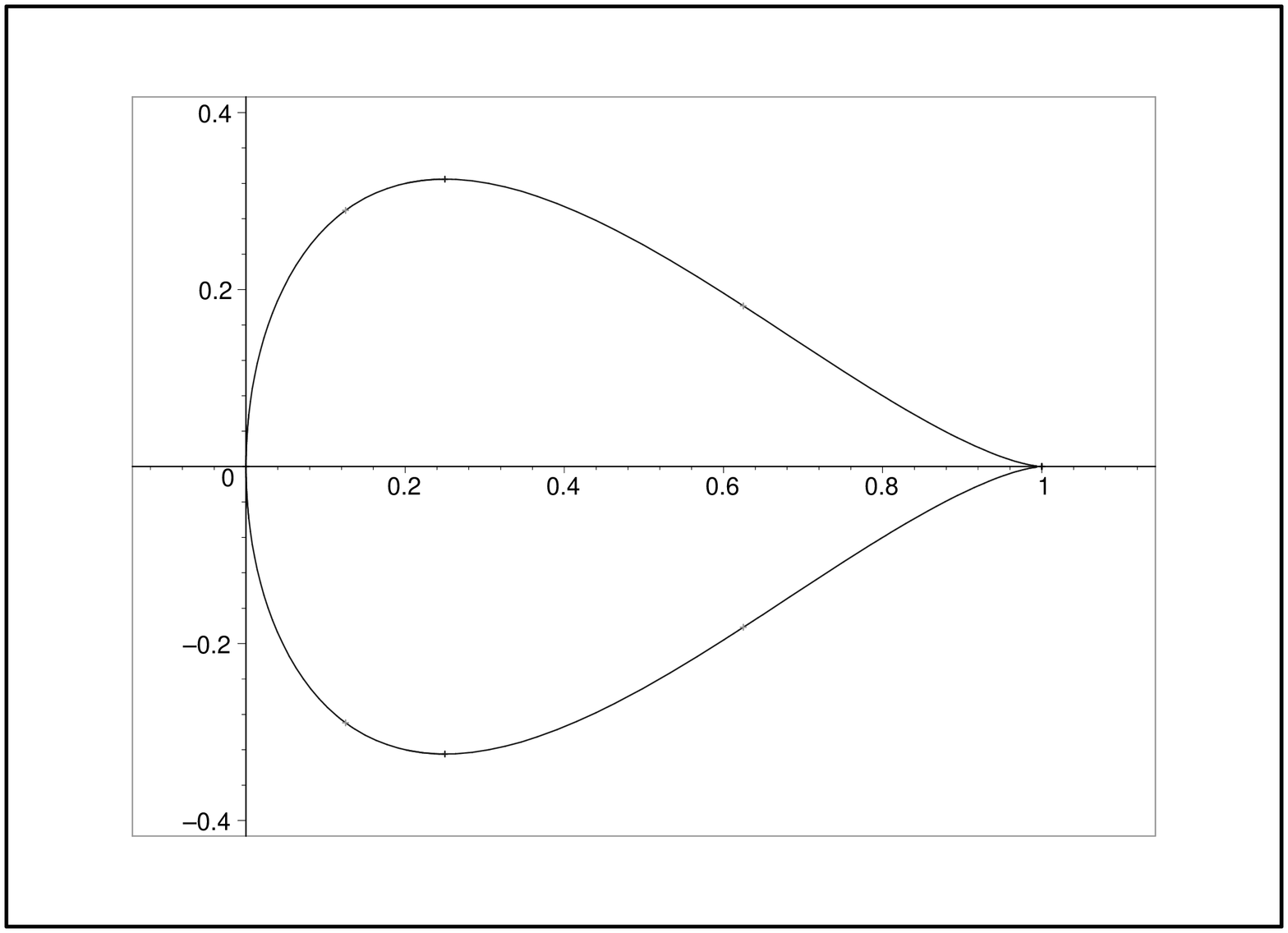}
\end{center}
\end{figure}
The curve (\ref{new_curve}) possesses the following property: its discriminant is a cube of a quartic $Q(x)$ 
in $x$ with two real roots. These roots correspond to the positions of
the finite cusps, as in Figure \ref{IS_levelsurface}. For any such  curve (\ref{sol_introduction}) the cross-ratio $\rho$ of the roots of $Q$ is an invariant with respect to transformations (\ref{point_tr}). It is possible to verify that for (\ref{new_curve}) we have $\rho=\varepsilon,$ where $\varepsilon=\exp{(\frac{i \pi}{3})}.$ To bring the curve 
(\ref{new_curve}) to a canonical form we first set 
the coefficient of $y^2$ to zero by a transformation of the 
form $y\rightarrow y+ c_1+c_2x+c_3x^2.$ After that the 
remaining M\"obius transformations can be used to 
set $Q(x)=x(x-1)(x-\rho).$ This,  up to scaling of $y$, yields
$$
y^3-3 y \Big(x-(1-\varepsilon) \Big)\Big(x-(1+\varepsilon) \Big)\Big(x-\frac{1}{3}(1+\varepsilon) \Big)\Big(x-(-1+\varepsilon) \Big)+
$$
$$
2 \Big(x-(1-i \varepsilon) \Big) \Big(x-(1+i \varepsilon) \Big)\Big(x-(i-i \varepsilon) \Big) \Big(x-(-i+i \varepsilon) \Big) \Big(x-(i+\varepsilon) \Big) \Big(x-(-i+\varepsilon) \Big)=0.
$$
The discriminant $Q(x)$ of this curve is proportional to $x^3 (x-1)^3 (x-\varepsilon)^3.$

Using this complex canonical form, we easily can find a real one. 
To do that we find four roots of the sextic in $x$ in the complex 
canonical form such that their cross-ratio is $-1$. 
Using the M\"obius transformations, we move these roots
to $\pm 1\pm i$. The result (up to scaling of $x$ and $y$) is given by 
\be
\label{canform}
y^3+3 (3 x^4-6 x^2-1) y+ 12 x (3 x^4+1)=0
\ee
which can be parametrised by
\[
x(t)=\frac{t(t^2-3)}{3(t^2+1)},\quad y(t)=-\frac{4t(t^4+3)}{3(t^2+1)^2}.
\]
 
To get the general solution  for the 7th order ODE (\ref{equation}) we 
apply the point transformations (\ref{point_tr}) to (\ref{canform}).  
The resulting degree six rational curve of the form 
(\ref{sol_introduction}) is given by
\begin{eqnarray}
\label{final_formula}
&& \left( c_{{4}}y+c_{{1}}+c_{{2}}x+c_{{3}}{x}^{2} \right) ^{3}\\
&&+3\,
 \left( 3\, \left( c_{{5}}x+c_{{6}} \right) ^{4}-6\, \left( c_{{5}}x+c
_{{6}} \right) ^{2} \left( 1-c_{{7}}x \right) ^{2}- \left( 1-c_{{7}}x
 \right) ^{4} \right)  \left( c_{{4}}y+c_{{1}}
+c_{{2}}x+c_{{3}}{x}^{2}\right)\nonumber\\
&& +12\, \left( c_{{5}}x+c_{{6}} \right)  \left( 3\, \left( c_{{5}}x+c_{{6}} \right) ^{4} \left( 1-c_{{7}}x \right) + \left( 1-c_{{7}}
x \right) ^{5} \right) =0\nonumber.
\end{eqnarray}
The formula  (\ref{sol_introduction}) is the general solution
as the generic initial data \[\{y(0), y'(0), \dots, y^{(6)}(0)\}\] can be 
chosen arbitrarily by choosing the coefficients $c_1, \dots, c_7$  in 
the solutions. There exist 
additional singular solutions corresponding to submanifolds of the initial
data manifold which can not be obtained by any choice of coefficients in 
(\ref{final_formula}). An example of such singular solution
with a co--dimension one initial data set is the curve (\ref{degree_six}).
The leading term in this curve is $y^2$ and this can not arise in
(\ref{final_formula}). This type of behaviour should not be confused with
the singular orbits of the symmetry group. The connected component of the 
group $Sp(4)$ is transitive on the solution space, and
thus (\ref{final_formula}) and (\ref{degree_six}) are related by a contact 
transformation. This is indeed how we found (\ref{final_formula}).

\subsection*{Conclusions}
We have found the general solution of the 7th order ODE (\ref{equation}). This equation is submaximal as its symmetry algebra is ten--dimensional, 
wheres the algebra of the trivial 7th order ODE is eleven--dimensional. 
The general solution to (\ref{equation}) is given by the seven
dimensional orbit of the point transformations (\ref{point_tr}) acting
on the canonical solution (\ref{canform}). There also exists a 
six--dimensional submanifold in the space of solutions given by degree four
curves ({\ref{degree_six}}).
For each contact transformations
generated by $H_8, H_9$ and $H_{10}$ in (\ref{10Dalgebra}) there exists a 
point transformation such that the composition of the two fixes  
(\ref{canform}). This gives a stabiliser $SL(2)$ of (\ref{canform})
and finally the solution space $Sp(4)/SL(2)$.

The analogous, submaximal
5th order ODE characterises conics in $\RP^2$. 
In the inhomogeneous coordinates $(x, y)$ 
the five parameter family of conics is
\[
y^2=c_1\;x^2+c_2\;xy+c_3\;y+ c_4\;x+c_5.
\]
Eliminating the parameters $(c_1, \dots, c_5)$ between this equation
and its fourth derivatives and substituting in the fifth derivative
yields the ODE
\[
9(y^{(2)})^2y^{(5)}-45y^{(2)}y^{(3)}y^{(4)}+40(y^{(3)})^3=0.
\]
This construction goes back to 
Halphen \cite{Halphen} who wrote the equation as $((y^{(2)})^{-2/3})^{(3)}=0$.
\vskip5pt
The seven--dimensional space of solutions $M$ to (\ref{equation})
carries a $GL(2, \R)$ structure in the sense of 
\cite{DT, Doubrov2}: $TM$ has a pointwise identification with a 
vector space of homogeneous degree six polynomials in two variables. 
The five W\"unschmann--Doubrov--Wilczynski invariants vanish on the ODE (\ref{equation}) which implies \cite{Doubrov1} that the
linearisation of (\ref{equation}) is equivalent to a trivial 
ODE $\delta y^{(7)}=0$. 
Moreover \cite{DGN10}, $M$ also admits a conformal structure of signature $(3, 4)$. 
The null vectors of this structure
correspond to the six order polynomials 
\[
a_1x^6+6a_2x^5+15a_3x^4+20a_4x^3+15a_5x^2+6a_6x+a_7
\]
with vanishing quadratic invariant \cite{Grace_Young}
\[
a_1a_7-6a_2a_6+15a_3a_5-10a_4^2.
\]
In \cite{DGN10} it was shown that the conformal structure associated to (\ref{equation}) contains a metric with weak holonomy $\tilde{G}_2$: there exists 
a three form $\phi$
on $M$ such that 
\[
d\phi=\Lambda *\phi, \quad d*\phi=0
\]
where $\Lambda$ is a constant, and $*$ is the Hodge operator. Taking an analytic continuation of this structure to the Riemannian signature yields the homology seven--sphere $M=SO(5)/SO(3)$
with its canonical weak $G_2$ structure originally constructed by 
Bryant \cite{Bryant1}.
\section*{Acknowledgements}
MD is grateful to  Robert Bryant, Boris Doubrov, Michal Godli\'nski and Peter Olver
for useful correspondence. VS was partially supported by the RFBR grant 08-01-00400. We thank the anonymous referee for perspective comments
on the manuscript. 

\end{document}